\documentclass{pspum-l}

\usepackage{palatino, mathpazo}
\usepackage{amsfonts, amssymb}
\usepackage{amsmath}
\usepackage[mathscr]{eucal}
\usepackage[active]{srcltx}
\usepackage[all]{xy}
\usepackage{color}
\usepackage{mathrsfs}

\newtheorem{theorem}{Theorem}[section]
\newtheorem{lemma}[theorem]{Lemma}
\newtheorem{proposition}[theorem]{Proposition}

\theoremstyle{remark}
\newtheorem{remark}[theorem]{Remark}
\newtheorem{example}[theorem]{Example}

\numberwithin{equation}{section}

\newcommand{\M}{\mathscr{M}}
\newcommand{\T}{\mathscr{T}}
\newcommand{\p}{\partial}

\title{Quantum Cohomology under Birational Maps and Transitions}
\author[Y.-P.~Lee]{Yuan-Pin~Lee}
\address{Y.-P.~Lee: Department of Mathematics, University of Utah,
Salt Lake City, Utah 84112-0090, U.S.A.}
\email{yplee@math.utah.edu}

\author[H.-W.~Lin]{Hui-Wen~Lin}
\address{H.-W.~Lin: Department of Mathematics and Taida
Institute for Mathematical Sciences (TIMS),
National Taiwan University, Taipei 10617, Taiwan}
\email{linhw@math.ntu.edu.tw}

\author[C.-L.~Wang]{Chin-Lung~Wang}
\address{C.-L.~Wang: Department of Mathematics and Taida
Institute for Mathematical Sciences (TIMS),
National Taiwan University, Taipei 10617, Taiwan}
\email{dragon@math.ntu.edu.tw}

\begin{document}

\maketitle

\begin{abstract}
This is an expanded version of the third author's lecture in String-Math 2015 at Sanya.
It summarizes some of our works in quantum cohomology. 

After reviewing the quantum Lefschetz and quantum Leray--Hirsch, 
we discuss their applications to the functoriality properties under special smooth flops, flips and blow-ups. 
Finally, for conifold transitions of Calabi--Yau 3-folds, formulations for small resolutions (blow-ups along Weil divisors) are sketched. 
\end{abstract}


\setcounter{section}{-1}

\section{Introduction}

\subsection{Classical aspects on algebraic geometry}

In the study of algebraic geometry, we usually encounter \emph{projective morphisms} $f: Y \longrightarrow X$, e.g.~blow-ups, bundle morphisms etc., and it is a basic question to study relations of geometric quantities under such morphisms.  

By its very definition of being projective, there are \emph{factorizations} of $f$ into compositions of simpler morphisms in the following form. There are vector bundles $\mathscr{E} \to X$ and associated projective bundles $\pi: P = P_X(\mathscr{E}) \to X$ such that $f = \pi \circ \iota$:
$$
\xymatrix{Y \ar[r]^>>>>>\iota \ar[rd]_f & P \ar[d]^\pi\\ & X}
$$
where $\iota: Y \hookrightarrow P$ is an imbedding. The choice of $P$ is by no means unique. In fact $\mathscr{E}$ can be taken to be a trivial bundle, say of rank $r$, and then $P = P^{r - 1} \times X$ is a product. However, a good choice of $P$ is usually important so that the induced imbedding $\iota$ has good structures. 

If $\iota(Y) \subset P$ is a complete intersection, namely that there is a split vector bundle $V = \bigoplus L_i \to P$ and a section $\sigma \in \Gamma(P, V)$ such that $ \iota(Y) = \sigma^{-1}(0)$ is its zero loci, then one develops 
$$
\mbox{Lefschetz Hyperplane Theorem}
$$
to study relations between $Y$ and $P$. In most cases one does not obtain complete intersection imbedding automatically. Nevertheless sometimes one may employ the technique of \emph{deformations to the normal cone} to reduce the problem under study to such a situation. The most famous one is the proof of \emph{Grothendieck--Riemann--Roch theorem} \cite{Fulton}. Historically it is the proof of GRR which lays the foundation of these techniques. Among many later applications, it is notable that \emph{deformations to the normal cone} is now a standard method to reduce a problem to its local models (e.g., the proof of \emph{invariance of complex elliptic genera} under $K$-equivalent birational maps \cite{Wang3}). 

For the (projective) bundle map $\pi$, one develops
$$
\mbox{Leray--Hirsch Theorem}
$$
to study relations between $P$ and $X$. The techniques in this package include direct images, spectral sequences etc., and usually can be applied to more general bundle maps such that the fiber manifold is well understood. A combination of results for $\iota$ and $\pi$ then gives the desired result for $f$.

\subsection{Quantum aspects}
In this article, following the above flow-chart, we survey the related developments on the quantum cohomology ring $QH(X)$, or equivalently the genus zero Gromov--Witten theory. More precisely we consider the Dubrovin (flat) connection $\nabla$ on $TH(X)$ and analyze its behavior under various maps including complete intersection imbedding and projective bundle maps. The essential mathematical tools are the corresponding \emph{Quantum Lefschetz Hyperplane Theorem} and the \emph{Quantum Leray--Hirsch Theorem}. 

The first version of quantum Lefschetz was proved around 1996 by Lian, Liu, and Yau \cite{LLY1} and by Givental \cite{Givental} independently (cf.~\cite{CoKa}). It is also known as the \emph{mirror theorem} since its major motivation and application is to prove the counting formula of rational curves on quintic Calabi--Yau 3-folds predicted by Candelas et.\ al.. In that setup $P$ could be a semi-Fano toric manifold and $V = \bigoplus L_i$ is a sum of convex line bundles such that $c_1(Y) \ge 0$ (semi-Fano or sub Calabi--Yau). There are several improvements of quantum Lefschetz afterwards (e.g.~\cite{ypL}). The most general version which allows the background manifold $P$ to be general and without the condition on $c_1(Y)$ was obtained by Coates and Givental \cite{CG}. We will review this latest theory in terms of Dubrovin connections and introduced the notion of \emph{Birkhoff factorizations} and \emph{generalized mirror transforms}.    

For quantum Leray--Hirsch, a version for $P = P_X(\mathscr{E})$ with $\mathscr{E} = \mathscr{O} \oplus L$ a rank two split bundle first appeared in the work of Maulik and Pandharipande \cite{MP}. A version on the more general case of toric bundle $P \to X$ build on a \emph{split vector bundle} $\mathscr{E} \to X$ was proved by Brown \cite{jB}. It had been formulated in the framework of Dubrovin connections and used to prove the \emph{invariance of quantum cohomology rings under ordinary flops} by us in 2011 \cite{LLWp1, LLWp2}. More recently, together with F.\ Qu, we proved a \emph{quantum splitting principle} and remove the splitting assumption in the quantum Leray--Hirsch \cite{LLQW}. For the applications to be discussed here, we will only focus on the Dubrovin connection in the split case. The ideas of \emph{n\"aive quantization basis} and \emph{admissible lift of Mori cone} are the essential ingredients to formulate the quantum Leray--Hirsch. 

By combining quantum Lefachetz and quantum Leray--Hirsch, we discuss several applications on birational maps including 
\begin{itemize}
\item[(i)] the analytic continuations, i.e.~invariance, of $QH$ under ordinary flops, 
\item[(ii)] smooth blow-ups along \emph{complete intersection} centers, as well as 
\item[(iii)] \emph{decomposition theorem} of projective local models of simple ordinary flips. 
\end{itemize}
The latter two applications are new, and both address the issue of the \emph{functoriality} of $QH$ under non $K$-equivalent transformations. 

Indeed for simple $(r, r')$ flips $f: X \dasharrow X'$ with $r > r'$, there is an orthogonal decomposition 
$$
H(X) = \T^{-1} H(X') \oplus K
$$ 
where $\T = [\Gamma_f]_*: H(X) \to H(X')$ is the map induced from the graph correspondence with kernel (vanishing cycles) $K \cong \Bbb C^{r - r'}$, and $\T^{-1} := [\Gamma_f^t]$ is the transpose correspondence. Under a suitable choice of quantum frame $\{\tilde T_i\}$ which deforms the classical cohomology basis $\{T_i\}$, we have a ring isomorphism (decomposition)
$$
QH(X) \cong \langle \tilde T_1, \ldots, \tilde T_{\dim H(X')}\rangle \times \Bbb C^{r - r'}
$$
such that $\langle \tilde T_1, \ldots, \tilde T_{\dim H(X')}\rangle \cong QH(X')$ as $\mathscr{D}^z$-modules, but not as rings. In this survey we illustrate this for the case of $(2, 1)$ flips.

Comparing with the current developments on string-math related topics, the subjects discussed here are more classical in flavor as there is essentially no higher genus theory nor modern mirror symmetry involved. Mirror symmetry phenomenon happens in \emph{the large complex/K\"ahler structure limit} while birational maps are essentially located at \emph{finite distance K\"ahler degenerations}. Also we treat only smooth varieties. However, in higher dimensional birational geometry, namely the \emph{minimal model theory}, it is indispensable to include singularities in the variety under consideration. 

It is clear that most of the subjects discussed here can be extended to orbifolds since orbifold Gromov--Witten theory is now well developed. However MMP requires more general singularities then orbifold ones and it is still a long way towards a useful \emph{quantum minimal model program}. Of course the smooth case is the first step, and it is our hope that further progress in the general case can be made in the near future.

\subsection{Towards a QMMP}
As the next step, we are led to consider \emph{birational maps up to complex deformations}, i.e.~\emph{transitions}. There are many technical issues from classic algebraic geometry on this regard. Nevertheless it has become clear in recent years that it is indispensable to allow certain transitions in the classification of higher dimensional varieties. 

For Calabi--Yau 3-folds the famous Ried's fantasy \cite{mR} on connecting CY with different topology through transitions is still one of the major research problems in this area. We give a very brief sketch on our recent work \cite{LLW3} in understanding the transition of quantum $\mathscr{A}$ model and $\mathscr{B}$ model in \emph{projective} conifold transitions $X \nearrow Y$ through a conifold degeneration $\mathfrak{X} \to \Delta$ from $X = \frak{X}_t$, $t \ne 0$, to $\bar X = \mathfrak{X}_0$, followed by a small resolution $Y \to \bar X$. In particular the notion of linked GW invariants with respect to a set of vanishing spheres $S_i \subset X$ is introduced which corresponds to the non-exceptional GW invariants on $Y$. For GW invariants supported on exceptional curves, a \emph{basic exact sequence} shows that there is a local transition between them and the Yukawa couplings of the vanishing periods. 

This article ends with discussions on issues of \emph{effective computations} on the quantum transitions in terms of blow-up formula of GW invariants when the blow-up is along certain non-complete intersection (Weil) divisors.

\section{Quantum cohomology}

In this article, unless stated otherwise, all varieties $X, Y, P, \ldots$ under considerations are assumed to be smooth and projective over $\Bbb C$. The cone of effective one cycles (Mori cone) in $X$ is denoted by $NE(X)$.

\subsection{Dubrovin connection}

A general reference is \cite{CoKa}. We fix a cohomology basis $T_i \in H = H(X) := H^*(X, \Bbb C)$, with dual basis $\{T^i\}$. A general element in cohomology is denoted by $t = \sum t^i T_i$. The genus zero GW theory is encoded by its generating function (pre-potential) $F(t) = \langle\!\langle \rangle\!\rangle$, where for $a_i \in H(X)$,
$$
\langle\!\langle a_1, \ldots, a_m \rangle\!\rangle := \sum_ {\beta \in NE(X)} \sum_{n = 0}^\infty \frac{q^\beta}{n!} \langle a_1, \ldots, a_m, t, \ldots, t \rangle_{g = 0, m + n, \beta}.
$$
The (formal) Novikov variables $q^\beta$'s are inserted to avoid the issue of convergence.  
It also keeps track on the natural grading arising from the virtual dimension (i.e., the \emph{conformal} structure on the \emph{Frobenius manifold} $H(X)$, though we do not make use of this language in this article). 
In general $t$ is also treated as a formal variable (cf.~\eqref{e:J} below for $t \in H^0 \oplus H^2$).

Denote by $F_{ijk} = \p^3_{ijk} F = \langle\!\langle T_i, T_j, T_k \rangle\!\rangle$ the 3-point generating functions, and set $A_{ij}^k := \sum_l F_{ijl}\,g^{lk}$. Then the big quantum product at $t$ is defined by 
$$
T_i *_t T_j = \sum_k A_{ij}^k(t) \, T_k. 
$$
The product is associative due to the WDVV equations. Equivalently it corresponds to the flatness of the Dubrovin connection on $T H \otimes \Bbb C [\![q^\bullet]\!]$:
\begin{equation} \label{e:Dubrovin}
\nabla := d - \frac{1}{z} A \equiv d - \frac{1}{z} \sum_i dt^i \otimes A_i.
\end{equation}
The special role played by the $z$ parameter shows that $\nabla$ is flat if and only if $dA = 0 = A \wedge A$, which is equivalent to WDVV.
 
\subsection{$J$ function and cyclic $\mathscr{D}^z$-modules} \label{s:DJ}

Write $t = t_0 + t_1 + t_2$ with $t_0 \in H^0$ and $t_1 \in H^2$. The generating function of all genus zero GW invariants with at most one descendent insertion is organized as
\begin{equation} \label{e:J}
\begin{split}
J(t, z^{-1}) &:= 1 + \frac{t}{z} + \sum_{\beta, n, i} \frac{q^{\beta}}{n!} T_i \left\langle \frac{T^i}{z(z - \psi)}, (t)^n \right\rangle_{\beta}\\
&= e^{\frac{t}{z}} + \sum_{\beta \neq 0, n, i} \frac{q^{\beta}}{n!} e^{\frac{t_0 + t_1}{z} + (t_1.\beta)} T_i \left\langle \frac{T^i}{z(z - \psi)}, (t_2)^n \right\rangle_{\beta}
\end{split}
\end{equation}
where the fundamental class axiom (string equation) and the divisor axiom are applied to get the second equality. The important role played by the $J$ function comes from the following quantum differential equation (QDE)
\begin{equation} \label{e:QDE}
z\p_i \, z\p_j\, J = \sum_k A_{ij}^k \, z\p_k\, J,
\end{equation}
which follows from the topological recursion relation (TRR). It implies that the quantum cohomology $QH(X)$ can be regarded as the cyclic $\mathscr{D}^z$-module $\mathscr{D}^z\,J$ with base (frame) given by
$$
z\p_i\, J \equiv e^{t/z}\,T_i \pmod {q^\bullet} = T_i + \ldots.
$$
It is clear that for $T_{1} = {\bf 1}$ being the fundamental class, $z\p_{1} \, J = J$.

The ring of differential operators $\mathscr{D}^z = \Bbb C[\![t, q^\bullet]\!]\{z, z\p_\bullet \}$ is defined so that $p = \sum_{\beta} q^\beta \,p_\beta \in \mathscr{D}^z$ implies that $p_\beta$ is a polynomial in $z$ and $z\p_\bullet$.

In practice it is sometimes easier to study the GW theory or the $J$ function on the small parameter space $H^0 \oplus H^2$. The expression of $J$ with $t_2 = 0$ is known as the small $J$ function. Using the divisorial reconstruction theorem in \cite{LP1}, the small $J$ function determines the sub-algebra of $QH(X)$ generated by $H^2(X)$.

\section{Review on quantum Lefschetz}

\subsection{Quantum Lefschetz for toric base and concavex bundle with $c_1 \ge 0$}

Let $P$ be a projective manifold, $L_i \to P$, $1 \le i \le r$ be convex line bundles, and $\sigma \in \Gamma(P, \bigoplus\nolimits_{i = 1}^r L_i)$ be a section such that $Y = \sigma^{-1}(0) \hookrightarrow P$ is a smooth submanifold. Given $QH(P)$, the problem is to compute $QH(Y)$. 

The method of localizations on stable map moduli spaces leads to the so called \emph{factorial trick} or \emph{hypergeometric modifications}. To state it, we start with the cohomology valued factorial 
$$
(L)_\beta := \prod_{m = 1}^{L.\beta} (L + mz)
$$
whenever the intersection number $L.\beta \ge 0$. Then we set
\begin{equation} \label{e:hg}
I^Y(t, z, z^{-1}) := \sum_{\beta \in NE(P)} q^\beta J^P_\beta(t, z^{-1}) \times \prod_{i = 1}^r (L_i)_\beta
\end{equation}
as an approximation of $J^Y$. 

When the ambient space $P$ is a semi-Fano toric manifold (i.e.~$c_1(P) \ge 0$), Lian--Liu--Yau and Givental used $\Bbb C^\times$ localizations to determine $J^P$. (This in turn determines $QH(P)$ since $H(P)$ is generated by divisors.) Over such a semi-Fano toric base, they further proved the   

\begin{theorem} [Mirror Theorem] \cite{LLY1, Givental} \label{t:mirror}
For $c_1(Y) \ge 0$, $t \in H^0 \oplus H^2$, we have
$$
(I^Y/I^Y_0) (t, z^{-1}) = J^Y(\tau, z^{-1})
$$ 
up to the mirror map $t \mapsto \tau(t)$ which matches $1/z$ coefficients on both sides.
\end{theorem}

Here $I^Y_0$ is the component of $z^0$ terms. Notice that the assumption $c_1(Y) \ge 0$ implies that $I^Y$ is still an expression in $z^{-1}$. 

In \cite{LLY1}, the line bundles $L_i$'s are also allowed to be concave. In that case their \emph{mirror principle} determines a certain type of twisted GW invariants.

\subsection{Quantum Lefachetz over general base and split bundles} \label{s:QL}

Without the condition that $c_1(Y) \ge 0$, the approximation $I^Y$ might contain terms with positive $z$ powers. But $J^Y$, by definition, contains only terms in powers of $z^{-1}$. Hence a more sophisticated transformation is needed in order to relate $I^Y $ to $J^Y$.

Coates and Givental considered the following situation. Let $P$ be a general projective manifold whose big quantum cohomology ring $QH(P)$ is given. Let $L_i \to P$, $1 \le i \le r$, be line bundles. They defined twisted GW invariants in this setup. When $L_i$'s are base-point free the twisted invariants are the GW invariants of the complete intersection sub-manifold $Y = \sigma^{-1}(0)$ for a generic section $\sigma \in \Gamma(P, \bigoplus_{i = 1}^r L_i)$. They proved

\begin{theorem} \cite{CG}
Given $J^P(t)$ in $t \in H(P)$, we have $I^Y \in \mathscr{D}^z J^Y$. More precisely, there exists a linear differential operator $b = b(\tau, z, q^\bullet, z\p_\bullet)$ which is polynomial in $z$ on each finite truncation of the Novikov variables $q^\bullet$ such that
$$
I^Y(t, z, z^{-1}) = b(\tau, z, q^\bullet, z\p_\bullet)\, J^Y(\tau, z^{-1}).
$$
Here $t \mapsto \tau(t) \in H(P)$ is a transformation determined by this property.
\end{theorem}

To get a better understanding of the statement, we notice that $b$ can be taken to be linear in $z\p_i$'s because of the QDE \eqref{e:QDE} on $J^Y$. In a similar fashion, the differentiation $z\p_i\, I^Y$ can also be represented by $\sum_j B_{ij}\,z\p_j J^Y$ for some formal functions $B_{ij}$'s which are polynomials in $z$ in any $\beta \in NE(P)$. This leads to the so called \emph{Birkhoff factorization} 
\begin{equation} \label{e:BF}
(z\vec\p\, I)(t, z, z^{-1}) = (z \vec \p\, J)(\tau, z^{-1}) B(\tau, z)
\end{equation}
of the square matrix $(z\vec\p\, I) = (z\p_1\,I, \ldots, z\p_{\dim H(P)}\,I)$.

Since $z\p_i\, I \equiv e^{t/z} T_i \equiv z\p_i \, J \pmod {q^\bullet}$, we have $B \equiv Id \pmod{q^\bullet}$ and this implies the isomorphism on $\mathscr{D}^z$ modules
$$
\mathscr{D}^z\,I (t) \cong \mathscr{D}^z \, J (\tau)
$$
up to a \emph{generalized mirror transform} $\tau(t)$ on $H(X)$. Now it is clear that
$$
p(t, z, q^\bullet, z\p_\bullet)\, I^Y(t, z, z^{-1}) = J^Y(\tau, z^{-1})
$$
for some linear operator $p$. In fact this operator plays the role to remove the $z^{\ge 0}$ terms in $I^Y$ and it can be effectively constructed by induction on $NE(P)$ (cf.~\cite[Theorem 1.10]{LLWp2} for a related construction). The map $t \mapsto \tau(t)$ is then determined by matching the $1/z$ coefficients on both sides.

We also notice that from \eqref{e:BF} the matrix $B^{-1}$ is the gauge transformation to bring the wrong frame $z\p_i\, I$'s back to the preferred frame $z\p_i\, J$ so that the connection matrix takes the form expected in \eqref{e:Dubrovin}. 


\section{Quantum Leray--Hirsch}

\subsection{Factorial trick for split projective (toric) bundles}

Let 
$$
\pi: P = P_X(V) \to X
$$ 
be a projective bundle. The classical Leray--Hirsch theorem computes the cohomology of the total space $P$ in terms of the base $X$ and the fibers. Let $h = c_1(\mathscr{O}_P(1))$, then 
$$
H(P) \cong \pi^* H(X)[h]/(f_V(h))
$$
where $f_V(h)$ is the Chern polynomial of the vector bundle $V \to X$.

It is natural to ask for a similar description on quantum cohomology. For this purpose we assume that $QH(X)$ is given, and $V = \bigoplus_{i = 1}^r L_i$ is a sum of line bundles. We seek for an analogous factorial trick as in the case of quantum Lefschetz. However the formulation must be different since now $QH(P)$ contains additional variables.

Let $\bar t \in H(X)$ be a general element from the base, $D = t^h h$ be the fiber divisor class with coordinate $t^h$, and we consider the mixed variable
$$
\hat t = \bar t + D.
$$
Then a hypergeometric modification of $J^X$ is defined by
\begin{equation} \label{e:HG}
I^P(\hat t, z, z^{-1}) = \sum_{\beta \in NE(P)} q^\beta J^X_{\pi_* \beta}(\bar t) \times e^{\frac{D}{z} + (D.\beta)}\prod_{i = 1}^r \frac{1}{(h + L_i)_\beta}. 
\end{equation}
Here the convention on factorial $\prod_{1}^s := \prod_{\infty}^{s}/\prod_{-\infty}^{0}$ is used so that $1/(L)_\beta$ makes sense even if $L.\beta < 0$. 

When $X = {\rm pt}$, this is the $I$ function of $P^{r - 1}$ coming from localizations. In general \eqref{e:HG} arises from fiber localization of a $\Bbb C^\times$-action, which exists by the split assumption on $V$. The formulation works for other fiber bundles as long as the fiber localization is well understood. In that case $D = \sum_{i = 1}^\rho t^i D_i$ for $D_1, \ldots, D_\rho$ being a basis of $H^2(P/X)$. The following is due to Brown:

\begin{theorem} \cite{jB} \label{t:brown}
Given $J^X(\bar t)$ in $\bar t \in H(X)$, we have $I^P \in \mathscr{D}^z J^P$. More precisely, there is a linear differential operator $b = \sum_{\beta \in NE(P)} q^\beta b_\beta$ with $\deg_z b_\beta < \infty$, and a graph $\hat t \mapsto \tau(\hat t): H(X) \oplus \Bbb Ch \to H(P)$, such that 
$$
I^P(\hat t, z, z^{-1}) = b(\tau, z, q^\bullet, z\p_\bullet) \, J^P(\tau, z^{-1}).
$$
\end{theorem}

\begin{remark}
The result was proved in \cite{jB} for split toric bundles, and stated in the language of \emph{Lagrangian cones}. We have presented it in an equivalent form to avoid introducing this machinery. 
\end{remark}

\subsection{The Dubrovin connection} \label{s:QLH}

Based on Theorem \ref{t:brown}, we had developed a method to compute the Dubrovin connection on the bundle space $P$ \cite{LLWp2}. It is roughly represented by following implication: 
$$
{\rm PF}^{P/X} + \nabla^X \Longrightarrow \nabla^P.
$$
To be precise, we need to introduce a system of equations controlling both the fiber directions and the base directions. 

For the fibers, we introduce the \emph{Picard--Fuchs ideal}. For the primitive fiber curve class $\ell \in NE(P/X)$, it is easily checked that $\Box_\ell I = 0$ where
\begin{equation} \label{e:PF-l}
\Box_\ell = \prod\nolimits_{i = 1}^r z\p_{h + L_i} - q^\ell e^{t^h}
\end{equation}
is the Picard--Fuchs operator. (Here $\p_L$ is the directional derivative in direction $L$.)
The Picard--Fuchs ideal is the left ideal of $\mathscr{D}^z$ generated by the Picard--Fuchs operators. 

For differentiations in the variables corresponding to the base, i.e.~$H(X)$, directions, we introduce \emph{the lifting of QDE} from $X$ to $P$. 
For each $\bar \beta \in NE(X)$, a lift of $\bar \beta$ is a curve class $\beta \in NE(P)$ such that $\pi_* \beta = \bar \beta$. Moreover, $\beta$ is called \emph{admissible} if
$$
-(h + L_i).\beta \ge 0, \quad \mbox{for all} \quad i = 1, \ldots, r.
$$
Admissible lift exists. In fact a minimal lift in effective classes is admissible.

Let $\bar \beta^* \in NE(P)$ be an admissible lift of $\bar\beta \in NE(X)$. We define
$$
D_{\bar\beta^*}(z) := \prod_{i = 1}^r \prod_{m = 0}^{-(h + L_i).\bar\beta^* - 1} (z\p_{h + L_i} - mz).
$$
Then it can be shown by direct computations that (cf.~\cite[Theorem 3.6]{LLWp2})
\begin{equation} \label{e:QDE-lift}
z\p_i \, z\p_j \, I = \sum\nolimits_{k, \bar\beta} q^{\bar \beta^*} e^{D.\bar\beta^*}\,\bar A_{ij,\,\bar\beta}^k(\bar t)\,D_{\beta^*}(z)\,z\p_k \, I .
\end{equation}

\begin{remark}
For toric bundles with $\rho(P/X) > 1$, the admissible lift still exists, though not unique, and the lifting of QDE is independent of the choices of $\bar\beta^*$ modulo the Picard--Fuchs ideal. A case of double projective bundle (hence $\rho(P/X) = 2$) will be discussed in Lemma \ref{l:PF-isom}.
\end{remark}

Let $\bar t = \sum \bar t^i \,\bar T_i \in H(X)$, and ${e = h^l\, \bar T_i} \in H(P)$. The n\"aive quantization of $e$ is defined to be the operator
$$
\hat e = \p^{ze} := (z\p_h)^l \,z\p_{\bar T_i} = (z\p_{t^h})^l\, z\p_{\bar t^i}.
$$
The idea for doing so is clear: since we do not have a variable corresponding to $e$, we simply use $l$-th derivatives in $t^h$ to approximate the directional derivative in $h^l$. On the base direction we keep the first order derivative in direction $\bar T_i$ since the variable $\bar t^i$ is available.

Using \eqref{e:PF-l} and \eqref{e:QDE-lift} we get the first order system on the frame $\p^{ze}\,I$ (with $e$ runs through a basis of $H(P)$) over the variables $t^a = \bar t^i, t^h$: 
$$
z\p_a \, (\p^{ze} \, I) = (\p^{ze} \, I)\, C_a(\hat t, z).
$$

As in \S \ref{s:QL}, we have the Birkhoff factorization matrix $B$ such that 
$$
(\p^{ze} \, I)(\hat t, z, z^{-1}) = (z \vec \p\, J)(\tau, z^{-1}) B(\tau, z).
$$
In fact $B^{-1}$ is the gauge transformation to remove $z^{\ge 0}$ in $C_a(\hat t, z)$. Furthermore, the map $\hat t \mapsto \tau(\hat t)$ is uniquely determined by matching the $1/z$ coefficients of the first column of $(\p^{ze} \, I)B^{-1}$ with $J$: 
$$
J(\tau, z^{-1}) = z\p_{\bf 1}\, J = p(\hat t, z, \p^{z\bullet})\,I(\hat t, z, z^{-1}).
$$
Set $z = 0$ in the gauge transformation we find $-(z\p_a B)B^{-1} \mapsto 0$ and
\begin{equation} \label{e:gauge}
B_0\, C_{a; 0}\, B_0^{-1} (\hat t)= \sum_{i = 1}^{\dim H(P)} A_i(\tau(\hat  t)) \frac{\p \tau^i}{\p t^a}(\hat t),
\end{equation}
where $B_0 = B(z = 0)$ and $C_{a; 0} = C_a(z = 0)$. 

Since $\tau \equiv \hat t \pmod {q^\bullet}$, by the Mori cone induction and divisorial reconstruction we may then determine all the Dubrovin connection matrices $A_i(t)$ from \eqref{e:gauge}. Of course the computations involved are necessarily complicated and very demanding. In applying these results special attention is paid to maintain the structural information. We will demonstrate on this through a few applications.

\section{Application I: Ordinary flops}

\subsection{The statement} Let $f: X \dasharrow X'$ be a $P^r$ flop. That is, there are two vector bundles $F, F' \to S$ of the same rank $r + 1$, such that the exceptional loci $Z \subset X$ has the following projective bundle structure: 
$$
{\rm Exc}\,f = Z = P_S(F) \mathop{\longrightarrow}\limits^{\bar\psi} S.
$$ 
Moreover, the normal bundle of $Z$ in $X$ is given by
$$
N = N_{Z/X} \cong \bar\psi^* F' \otimes \mathscr{O}_Z(-1).
$$

It was shown in \cite{LLW} that the graph correspondence $\T = [\Gamma_f]_*$ induces an isomorphism on cohomology spaces $H(X) \cong H(X')$, but it does not preserve the effectivity of one cycles. Indeed if $\ell$ (resp.~$\ell'$) is the class of primitive extremal rational curve in $X$ (resp.~$X'$) then 
$$
\T \ell = -\ell'.
$$
Moreover, $\T$ preserves the Poincar\'e pairing, but not the product structure. It turns out that the topological defects are corrected by the extremal ray GW invariants and the following is true:

\begin{theorem} \cite{LLW, LLWp1, LLWp2, LLQW} \label{t:QH-isom}
The graph correspondence $\T$ induces isomorphism of big quantum cohomology rings $QH(X) \cong QH(X')$ under the analytic continuations induced from $q^\beta \mapsto q^{\T \beta}$.
\end{theorem}

Here is a brief history on this problem.
For $\dim X = 3$, the multiple cover formula for $\mathscr{O}_{P^1}(-1)^2 \to P^1$ gives the quantum corrections of $(\T D)^3 - D^3$ (cf.~Witten \cite{Witten}). The global case was treated by Li--Ruan \cite{LiRu} around 2000. They proved a degeneration formula of GW invariants in the \emph{symplectic category}, and used it to show that in fact no degeneration might occur in the threefold case. The problem was then reduced to the case of  extremal rays which had already been solved. 

In higher dimensions the statement was conjectured to hold for general \emph{birational $K$-equivalent manifolds} by Ruan and Wang (cf.~\cite{Wang2}). The case of \emph{simple} ordinary flops, namely $S = {\rm pt}$, was solved in 2006 \cite{LLW}. We worked in the algebraic category and reduced the problem to local models by way of the deformations to the normal cone and the degeneration formula of Li \cite{Li}. In this case non-trivial degenerations do arise and the relations between relative GW invariants, descendent invariants, and absolute GW invariants are carefully studied through degenerations (inspired by a method of Maulik and Pandharipande in \cite{MP}). For local models, $X = P_{P^r}(\mathscr{O}(-1)^{r + 1} \oplus \mathscr{O})$ is a semi-Fano toric variety whose GW theory is well studied (cf.~Theorem \ref{t:mirror}). In fact $I^X = J^X$ on small parameters, and the analytic continuation can be solved. The result was further extended to the higher genus GW theory in \cite{ILLW} by studying \emph{ancestor invariants} and \emph{quantization}.

For general base $S$ with split bundles $F$, $F'$, the analytic continuation was later solved in 2011. Indeed, the problem was reduced to the local models in \cite{LLWp1}, and the case of local models was solved through the quantum Leary--Hirsch theorem in \cite{LLWp2}. More recently, a \emph{quantum splitting principle} was proved in a joint work with Qu in \cite{LLQW}, which reduced the problem for general vector bundles $F, F'$ to the case of split bundles, hence proved Theorem \ref{t:QH-isom} completely.

\subsection{A sketch of proof} \label{s:sketch-pf}

Now we sketch how the Quantum Leray--Hirsch is applied to solve the case of local split flops. 

The flop is achieved by first blowing up $Z \subset X$ to get $Y = {\rm Bl}_Z X \to X$ and then contracting the exceptional divisor $E = P_S(F) \times_S P_S(F') \subset Y$ in another fiber direction to get $Z' \subset X'$. We have $\bar\psi': Z' = P_S(F') \to S$ being a projective bundle and $N' = N_{Z'/X'} = \bar\psi'^* F \otimes \mathscr{O}_{Z'}(-1)$:
$$
\xymatrix{
X= P_Z(N \oplus \mathscr{O}) \ar@{.>}[rr]^f \ar[rd]_{p} & & X' = P_{Z'}(N' \oplus \mathscr{O}) \ar[ld]^{p'} \\ & S},
$$
where $p = \bar\psi \circ \pi$ in $X \mathop{\longrightarrow}\limits^\pi Z \mathop{\longrightarrow}\limits^{\bar\psi} S$ and similarly $p' = \bar\psi' \circ \pi'$.
As a double projective bundle, we have $NE(X/S) = \langle\ell, \gamma \rangle$, where $\ell$ (resp.~$\gamma$) is the $\bar\psi$ (resp.~$\pi$) fiber line classes. Let $\xi = \mathscr{O}_X(1)$ and $h = \mathscr{O}_Z(1)$, and $D = t^h h + t^\xi \xi \in H^2(X/S)$ a general fiber divisor. Then
$$
H(X) = p^*H(S) [h, \xi]/(f_F, f_{N \oplus \mathscr{O}}),
$$ 
where $f_V$ is the Chern polynomial of a bundle $V$. 

When $F = \bigoplus_{i = 1}^r L_i$, $F' = \bigoplus_{i = 1}^r L'_i$ are split bundles, we have 
\begin{equation*}
\begin{split}
f_F(h) &= \prod (h + L_i), \\
f_{N \oplus \mathscr{O}}(h, \xi) &= \xi\prod (\xi - h + L'_i).
\end{split}
\end{equation*}
By symmetry we have similar formulae on the $X'$ side. However, it is not compatible with that on $X$: since $\T h =\xi' - h'$, $\T \xi = \xi'$, it is easy to see that the cup product structure is not preserved under $\T$. (See \cite[Theorem 1.8]{LLWp1} for the explicit computations on the topological defects.)

Now comes the key point: to remedy the topological defects we replace the cohomology class by its ``quantized'' version, namely we consider differential operators instead. This gives rise to the Picard--Fuchs operators
\begin{equation*}
\begin{split}
\Box_\ell &= \prod z\,\p_{h + L_i} - q^\ell e^{t^h} \prod z\,\p_{\xi - h + L'_i}, \\
\Box_\gamma &= z\,\p_\xi \prod z\,\p_{\xi - h + L'_i} - q^\gamma e^{t^\xi}
\end{split}
\end{equation*}
which are regarded as the ``quantized version'' of the Chern polynomials. Similarly we have on the $X'$ side
\begin{equation*}
\begin{split}
\Box_{\ell'} &= \prod z\,\p_{h' + L'_i} - q^{\ell'} e^{t^{h'}} \prod z\,\p_{\xi' - h' + L_i}, \\
\Box_{\gamma'} &= z\,\p_{\xi'} \prod z\,\p_{\xi' - h' + L_i} - q^{\gamma'} e^{t^{\xi'}}.
\end{split}
\end{equation*}
The coordinates are related by requiring $t^{h'} h' + t^{\xi'} \xi' = \T (t^h h + t^\xi \xi) = t^h (\xi' - h') + t^\xi \xi'$. That is, $t^{h'} = -t^h$ and $t^{\xi'} = t^\xi + t^h$.

Now it is a simple exercise to check that

\begin{lemma} \label{l:PF-isom}
$\T$ induces an isomorphism on Picard--Fuchs ideals 
$$
\T\langle \Box_\ell, \Box_\gamma \rangle \cong \langle \Box_{\ell'}, \Box_{\gamma'} \rangle.
$$  
\end{lemma}

\begin{remark}
Lemma \ref{l:PF-isom} can be extended to split toric bundle flops.
\end{remark}

For the base directions, the lift of QDE in \eqref{e:QDE-lift} is independent of the choice of admissible lift $\bar\beta^*$ modulo the Picard--Fuchs ideal. The admissible condition for $\beta$ in this case is given by $-\beta.(h + L_i) \ge 0$, $-\beta.(\xi - h + L'_i) \ge 0$ and $-\beta.\xi \ge 0$. It is readily seen that $\beta$ is admissible in $X$ if and only if $\T \beta$ is admissible in $X'$. This implies that the lifting of QDE from $S$ to $X$ and the one from $S$ to $X'$ are indeed equivalent under $\T$ modulo the Picard--Fuchs ideal. Thus, the quantum Leray--Hirsch theorem implies that $X$ and $X'$ have compatible first order PDE systems up to analytic continuations.

To achieve $\T: QH(X) \cong QH(X')$, we still need to show that the Birkhoff factorization $B$ and the generalized mirror map $\tau(\hat t)$, as appeared in \eqref{e:gauge}, are compatible on both sides. We refer the details to \cite[\S 3.3]{LLWp2}.

\section{Application II: Blow-ups along complete intersection centers}

Let $L_i = \mathscr{O}_X(D_i)$, $1 \le i \le r$, and $Z = D_1 \cap \cdots \cap D_r$ be a smooth complete intersection subvariety of $X$ with codimenison $r$. Let $\mathscr{E} = \bigoplus_{i = 1}^r L_i$ with a given section $s = (s_i)$ such that $D_i = (s_i)$, $Z = s^{-1}(0)$. Consider the blow-up $\phi:Y \to X$ along $Z$:
\begin{equation*} 
\xymatrix{E \ar[d] \ar@{^{(}->}[r] & Y = {\rm Bl}_Z X \ar[d]^\phi  \\ Z \ar@{^{(}->}[r] & X}.
\end{equation*}
Given $QH(X)$ and $(\mathscr{E}, s)$, the problem is to determine $QH(Y)$.

By construction, we have a surjective morphism $\mathscr{E}^* \twoheadrightarrow \mathscr{I}_Z$. This leads to the imbedding
$$
\iota: Y := {\rm Proj}_X \bigoplus_{d = 0}^\infty \mathscr{I}_Z^d \hookrightarrow {\rm Proj}_X \,{\rm Sym}\, \mathscr{E}^* = P_X(\mathscr{E}).
$$ 
Let $\pi: P := P_X(\mathscr{E}) \to X$ be the bundle map with associated Euler sequence $0 \to S \to \pi^*\mathscr{E} \to Q \to 0$ over $P$. It is shown in \cite{LLW4} that there is a canonical section $\sigma \in \Gamma(P, Q)$ of the universal quotient bundle such that $Y = \sigma^{-1}(0) \subset P$. We emphasize that this is not true if $Z$ is not a complete intersection. The situation is summarized in the following diagram:
\begin{equation} \label{e:Euler}
\xymatrix{0 \ar[r] & S \ar[rd] \ar[r] & \pi^*\mathscr{E} \ar[d] \ar[r] & Q \ar[ld] \ar[r] & 0\\ & Y \ar[rd]_\phi \ar@{^{(}->}[r]^>>>>>>\iota & P \ar[d]^\pi \ar@/_/[ru]_>>>>>\sigma \\ & & X}.
\end{equation}
Let $\eta = c_1(\mathscr{O}_P(1))$. It follows that $S = \mathscr{O}_P(-\eta)$ and $-\eta|_Y = E$ (cf.~\cite{Fulton}).

Notice that the quotient bundle $Q$ is in general not a split bundle over $P$ and the quantum Lefschetz can not be applied directly. However, a suitable extension of it to short exact sequences allows us to apply it to the Euler sequence \eqref{e:Euler} where both $S$ and $\pi^*\mathscr{E}$ are split bundles.

With the above understood, the (extended) quantum Lefschetz together with the quantum Leray--Hirsch lead to the following factorial trick on the $\beta$ component of $I^Y$:
\begin{equation} \label{e:QL+QLH}
\begin{split}
I^Y_\beta &= J^P_\beta \frac{\prod_{i = 1}^r (D_i)_\beta}{(-\eta)_\beta} \\ 
&\sim J^X_{\pi_* \beta} \,e^{\frac{t^\eta\,\eta}{z} + t^\eta(\eta.\beta)}\,\frac{\prod_{i = 1}^r (D_i)_\beta}{(-\eta)_\beta \prod_{i = 1}^r (\eta + D_i)_\beta}.
\end{split}
\end{equation}
Here $t^\eta$ denotes a suitable dual coordinate of $\eta.$

\begin{theorem} \cite{LLW4}
For the smooth blow-up $\phi: Y \to X$ along a complete intersection center $Z = \bigcap_{i = 1}^r D_i$, the relative $I$ factor is given by
\begin{equation*} 
I^{Y/X}_\beta = e^{s(\frac{E}{z} + E.\beta)} \left(\prod_{i = 1}^r\frac{(D_i)_\beta}{(D_i - E)_\beta (E)_\beta} \right) (E)_\beta^{r - 1},
\end{equation*}
where $E \subset Y$ is the exceptional divisor and $s$ is a suitable dual coordinate.
\end{theorem}

\begin{remark} \label{r:L-system}
The formula suggests nice structures of the relative factor. Indeed, $K_Y = \phi^* K_X + (r - 1) E$ and $(E)^{r - 1}_\beta$ is responsible for the Jacobian of $\phi$. For each $i$, $(D_i)_\beta/((D_i - E)_\beta (E)_\beta)$ is a \emph{$K$-trivial factor} which describes the decomposition of the linear system $|D_i - Z|$ into moving part and the fixed part (on $Y$). We expect that this intrinsic formulation will be useful for general blow-ups. 
\end{remark}

As in quantum Leray--Hirsch described in \S \ref{s:QLH}, to get the Dubrovin connection on $Y$ (or $QH(Y)$) we proceed by (1) choosing the corresponding n\"aive quantization basis (2) determining the Picard--Fuchs ideal on fibers (3) finding the lifting of QDE on the base $X$ to $Y$. Then $(1) + (2) + (3)$ determines $\mathscr{D}^z J^Y$, and hence $QH(Y)$. The details will appear in \cite{LLW4}. 


\section{Application III: Simple flips}

Let $r, r' \in \Bbb N$ and $r > r'$. In the definition of ordinary flops, if the underlying vector bundles $F$ and $F'$ have different rank $r + 1$ and $r' + 1$ respectively then in exactly the same construction as \S \ref{s:sketch-pf} we get ordinary $(r, r')$ flips. The effect on quantum cohomology under flips are discussed in \cite{LLW5}. In contrast to analytic continuations in the flops case, the situation for flips is more subtle and complex and new phenomena appear. We give a sketch in the simplest case of local models of simple $(2, 1)$ flips. 

\subsection{$H(X)$ vs $H(X')$ and the Picard--Fuchs systems}

The local model of $(2, 1)$ flips has the following data: $Z = P^2$, $Z' = P^1$, and
$$
f: X = P_{P^2}(\mathscr{O}(-1)^2 \oplus \mathscr{O}) \dasharrow X' = P_{P^1}(\mathscr{O}(-1)^3 \oplus \mathscr{O}).
$$
The cohomology rings are given by 
\begin{equation*}
\begin{split}
H(X) &= \Bbb C[h, \xi]/(h^3, \xi(\xi - h)^2), \\
H(X') &= \Bbb C[h', \xi']/(h'^2, \xi'(\xi' - h')^3),
\end{split}
\end{equation*}
where $\dim H(X) = 9$ and $\dim H(X') = 8$.

The graph correspondence $\T = [\Gamma_f]_*$ induces a short exact sequence
$$
0 \longrightarrow K \longrightarrow H(X) \mathop{\longrightarrow}^\T H(X') \longrightarrow 0
$$
where $K = \ker \T = \Bbb C {\bf k}_1$ with ${\bf k}_1 = (\xi - h)^2 = [Z]$. 

The transpose correspondence $\T^{-1} := [\Gamma_f^t]_*$ preserves the Poincar\'e pairing and induces an imbedding  $\T^{-1}: H(X') \hookrightarrow H(X)$ (indeed, of motives) which leads to an orthogonal decomposition (cf.~\cite[\S 2.3]{LLW})
$$
H(X) = \T^{-1} H(X') \mathop{\oplus}^\perp K.
$$
However, $\T^{-1}$ does not preserves the cup product. In fact $K^\perp$ is not closed under cup product. As in the case of flops we have curve classes $\ell, \gamma$ in $X$ and $\ell', \gamma'$ in $X'$. They are related by $\T \ell = -\ell'$ and $\T \gamma = \ell' + \gamma'$. Also $\T h = \xi' - h'$ and $\T \xi = \xi'$.

The divisor variable takes the form $D = t^h h + t^\xi \xi$. To simplify notations in our discussion, we will use variables 
$$
q_1 = q^\ell e^{t^h}, \quad q_2 = q^\gamma e^{t^\xi}; \qquad q_1' = q^{\ell'} e^{-t^h}, \quad q_2' = q^{\gamma'} e^{t^h + t^\xi}.
$$

From the computational point of view of quantum cohomology, $X'$ is bad since $c_1(X') = -h' + 4\xi'$ which contains both $K$-positive and $K$-negative directions. In other words, the expression $I^X$ is complicated and contains positive $z$ powers. Its Picard--Fuchs equations are given by
$$
\Box_{\ell'} = (z\p_{h'})^2 - q_1' (z\p_{\xi' - h'})^3, \qquad \Box_{\gamma'} = z\p_{\xi'}(z\p_{\xi' - h'})^3 - q_2'.
$$
The first equation $\Box_{\ell'} I^{X'} = 0$ shows that in order to reduce $(z\p_{h'}^2)I^{X'}$ we will receive derivatives of even higher power. It does give the correct reduction algorithm in the Mori cone topology since there is also a $q_1'$ multiplied. It is difficult to compute $\nabla^{X'}$ or to get any structure of it from this approach. 

On the other hand, $X$ is toric Fano with $c_1(X) = h + 3\xi$. The computation of $QH(X)$ is in principle easy since $I^X = J^X$ along small parameters. Assume that this has been done, then a natural question is 
$$
\mbox{Can we get $QH(X')$ from $QH(X)$ in a canonical manner?}
$$

Now we restrict ourselves to the small parameters $t \in H^0 \oplus H^2$ so that we can work with variables $q_1, q_2$ and $q_1', q_2'$ directly. The Picard--Fuchs equations on $J = J^X$ can be easily determined to be
$$
\Box_\ell = (z\p_h)^3 - q_1 (z\p_{\xi - h})^2, \qquad \Box_\gamma = z\p_\xi (z\p_{\xi - h})^2 - q_2.
$$
It is closely related to the one for $X'$: 

\begin{lemma} \label{l:PF-isom'}
Along the partially compactified two dimensional K\"ahler moduli $\mathscr{K} := \{(q_1, q_2)\} \bigcup \{(q_1' = 1/q_1, q_2' = q_1 q_2)\} \cong \mathscr{O}_{P^1}(1)$, the Hopf--Mobius stripe, we have
$$
\T: \langle \Box_\ell, \Box_\gamma \rangle \cong \langle \Box_{\ell'}, \Box_{\gamma'} \rangle
$$
outside the divisors $D_0 = \{q_1 = 0\}$ and $D_\infty = \{q_1' = 0\}$.
\end{lemma}


\subsection{Exact formula for $\nabla^X$ \cite{LLW5}}

The following frame (recall that $I = J$)
\begin{equation*} 
\begin{split}
v_1 &= \hat {\bf 1}J = J, \\
v_2 &= \hat hJ, \quad v_3 = (\hat \xi - \hat h)J, \\
v_4 &= \hat h^2J - (\hat \xi - \hat h)^2J, \quad v_5 = \hat h(\hat \xi - \hat h)J + (\hat \xi - \hat h)^2J, \\
v_6 &= \hat h^3J - \hat h(\hat \xi - \hat h)^2J,\quad v_7 = \hat h^2(\hat \xi - \hat h)J + \hat h(\hat \xi - \hat h)^2J, \\
v_8 &= \hat h^3 (\hat \xi - \hat h)J + \hat h^2(\hat \xi - \hat h)^2J, \\
v_9 &= \hat {\bf k}_1 J = (\hat \xi - \hat h)^2J,
\end{split}
\end{equation*}
respects $H(X) = \T^{-1}H(X') \oplus^\perp K$ when modulo $q_1, q_2$. They are precisely 
$$
z\p_i\,J \quad \mbox{at $t \in H^0 \oplus H^2$}, \quad 1 \le i \le 9,
$$ 
and we get the Dubrovin connection matrices 
\begin{equation*}
\begin{split}
A_1 = h*_{small} &= 
\begin{bmatrix}
&&&&& q_1 q_2 \\
1 \\
&&&&&&& q_1 q_2 \\
& 1 \\
&& 1\\
&&& 1 &&&&& -1 \\
&&&& 1 \\
&&&&& -1 & 1 \\
& 1 & -1 &&&&&& q_1
\end{bmatrix},\\
A_2 = \xi *_{small} &= 
\begin{bmatrix}
&&& -q_2 & q_2 & q_1 q_2 &&& q_2\\
1 &&&&& -q_2 & q_2\\
1 &&&&&&& q_1 q_2 \\
& 1 &&&&&& q_2\\
& 1 & 1\\
&&& 1 &&&&& \\
&&& 1 & 1 \\
&&&&&& 1 \\
&&&&&&& q_2 &
\end{bmatrix}.
\end{split}
\end{equation*} 

The crucial observation is that there is a single appearance of $q_1$ in $A_{1, 9}^9$. This shows that the system has irregular singularity along $D_\infty = (q_1 = \infty)$ in the $K$ direction. Let $x = q_1'$, $y = q_2'$. After a constant change of basis from $v_i$'s to $w_i$'s such that the Poincar\'e pairing 
$$
(w_i, w_j) = \delta_{9, i + j}, \qquad 1 \le i, j \le 8
$$
and $w_9 := v_9$ with $(w_9, w_i) = \delta_{9, i}$, the fundamental solution matrix $S$ satisfies  
\begin{equation*}
z(x\,\p_x) S = 
\begin{bmatrix}
&&& -\frac{1}{2} xy & xy &&&& xy \\
&&&&& -\frac{1}{2} xy & xy \\
1 &&&&& \frac{1}{4} xy & -\frac{1}{2} xy\\
&&&&&&& xy \\
& 1 &&&&&& -\frac{1}{2} xy\\
&&&&&&&& 1 \\
&&& 1 &&&&& -\frac{1}{2} \\
&&&&& 1 \\
& -\frac{1}{2} & 1 &&&&& xy & -1/x
\end{bmatrix} S,
\end{equation*} 
which is irregular in the $K$-block, i.e.~the (9, 9) entry, of Poincar\'e rank one.

\subsection{Block diagonalization}

The classical theory of ODE and the flatness of $\nabla^X$ imply that there exists a unique formal gauge transformation $S = P Z$:
\begin{equation}
P(x, y, z) = 
\begin{bmatrix} \label{e:P-gauge}
1 &&& g_1 \\
& \ddots && \vdots \\
&& 1 & g_8 \\
f_1 & \cdots & f_8 & 1
\end{bmatrix},
\end{equation}
such that 
$$
z (x\,\p_x) Z = B_1\, Z, \qquad z (y \,\p_y) Z = B_2\, Z,
$$ 
with $B_1$, $B_2$ being block-diagonalized. Moreover, each $f_i (x, y, z) = - g_{9 - i}(x, y, -z)$ is a formal series expansion of certain special function. The claim is that we may relate the first $8 \times 8$ blocks of $B_1(x, y, z)$ and $B_2(x, y, z)$ with the Dubrovin connection $\nabla^{X'}$. 

Under the new $z$-dependent frame $\tilde w_1, \ldots, \tilde w_8, \widetilde {\bf k}_1 J$ from \eqref{e:P-gauge}, namely
\begin{equation} \label{e:P-frame}
\tilde w_i = w_i + f_i \,{\bf \hat k}_1 J, \qquad {\bf \widetilde k}_1 J = {\bf \hat k}_1 J + \sum_{i = 1}^8 g_i\, w_i,
\end{equation}
we have (a special case of \cite{LLW5} for $(r, r') = (2, 1)$):

\begin{theorem}
For a simple $(2, 1)$ flip $f: X \dasharrow X'$, under the frame \eqref{e:P-frame} at $z = 0$, we have a ring isomorphism 
$$
QH(X) \cong \langle \tilde w_1(0), \ldots, \tilde w_8(0)\rangle \times \Bbb C.
$$
Moreover,
$\langle \tilde w_1(0), \ldots, \tilde w_8(0)\rangle \cong QH(X')$ as $\mathscr{D}^z$-modules, but not as rings.
\end{theorem}

The proof is based on Lemma \ref{l:PF-isom'} and we refer to \cite{LLW5} for the details.

\section{Conifold transiitons of CY 3-folds}

\subsection{Relations on vanishing $A$ and $B$ cycles} 

A Calabi--Yau variety is a $\Bbb Q$-Gorenstein variety with $K \sim 0$ and $h^1(\mathscr{O}) = 0$. 

Let $X \nearrow Y$ be a \emph{projective} conifold transition of Calabi--Yau 3-folds $X$, $Y$ through a singular Calabi--Yau variety $\bar X$ with $k$ ordinary double points (ODPs) $p_1, \ldots, p_k \in \bar X$. During the complex degeneration $\pi: \mathfrak{X} \to \Delta$ with $\mathfrak{X}_0 = \bar X$,  there are $k$ vanishing 3-spheres $S_1, \ldots, S_k$ with $N_{S_i/X} = T^* S^3$. During the K\"ahler degeneration (small contraction) $\psi: Y \to \bar X$, there are $k$ vanishing 2-spheres (exceptional curves) $C_1, \ldots, C_k$ with $N_{C_i/Y} = \mathscr{O}_{P^1}(-1)^{\oplus 2}$:
$$
\xymatrix{& {C_i \subset Y \ } \ar[d]^\psi \\ S_i \subset X \ar@{~>}[r]^\pi & p_i \in \bar X \ .}
$$

Let $\mu := h^{2,1}(X) - h^{2, 1}(Y) > 0$ be the lose of complex moduli and $\rho := h^{1, 1}(Y) - h^{1, 1}(X) > 0$ be the gain of K\"ahler moduli. From 
$$
\chi(X) - k\chi(S^3) = \chi(Y) - k \chi(S^2),
$$ 
we get the following well-known elementary relation 
$$
\mu + \rho = k.
$$ 
This implies that the $\psi$-exceptional curve classes $[C_i] \in NE(Y/\bar X)$ admit $\mu$ independent relations, and the $\pi$ vanishing cycles $[S_i] \in V \hookrightarrow H_3(X) \to H_3(\bar X)$ admit $\rho$ independent relations. (The vanishing cycle space $V$ has $\dim V = \mu$.) Let $A$, $B$ be the corresponding relation matrices:
\begin{equation*}
\begin{split}
A = (a_{ij}) \in M_{k \times \mu}, \qquad &\sum\nolimits_{i = 1}^k a_{ij} [C_i] = 0, \\
B = (b_{ij}) \in M_{k \times \rho}, \qquad &\sum\nolimits_{i = 1}^k b_{ij} [S_i] = 0.
\end{split}
\end{equation*} 

\begin{theorem} [Basic exact sequence] \cite[Theorem 1.14]{LLW3} \label{t:bes}
The Hodge realization of $\mu + \rho = k$ is represented by an exact sequence 
$$
0 \to H^2(Y)/H^2(X) \mathop{\longrightarrow}^B \Bbb C^k \mathop{\longrightarrow}^{A^t} V \to 0
$$
of weight two Hodge structures.
\end{theorem}

Indeed $V \cong H^{1,1}_\infty H^3(X)$ in the limiting Hodge diamond for $\pi$:
\begin{equation*}
\xymatrix{&&H^{2, 2}_\infty H^3 \ar[dd]^N_\sim\\
	\Bbb C \cong H^{3,0}_\infty H^3 &H^{2,1}_\infty H^3 \ar@{-}[ru]&&H^{1,2}_\infty H^3 \ar@{-}[ld]&H^{0,3}_\infty H^3 \\
	&&H^{1, 1}_\infty H^3}  
\end{equation*}
and the \emph{invariant subsystem} is $Gr^W_3 H^3(X) \cong H^3(Y)$.

\subsection{Local quantum transition} 

By the Bogomolov--Tian--Todorov theorem and its extension to Calabi--Yau conifolds by Ran and Kawamata, the moduli spaces $\M_Y$ and $\M_{\bar X}$ are smooth of dimension $h^{2, 1}(Y)$ and $h^{2, 1}(X)$ respectively. Also the contraction $\psi: Y \to \bar X$ deforms in projective families. This then identifies $\M_Y$ as a codimenison $\mu$ boundary strata in $\M_{\bar X}$ and \emph{locally} near $[\bar X] \in \M_{\bar X}$ we have $\M_{\bar X} \cong \Delta^\mu \times \M_Y$.

We represent $V = \Bbb C\langle \Gamma_1, \ldots, \Gamma_\mu\rangle$ in terms of a basis $\Gamma_j$'s. It was shown in \cite[Proposition 3.15]{LLW3} that the $\alpha$-periods 
$$
r_j = \int_{\Gamma_j} \Omega, \qquad 1 \le j \le \mu
$$ 
form the \emph{degeneration coordinates} around $[\bar X] \in \M_{\bar X} \cong \Delta^\mu \times \M_Y$.

In order to describe the discriminant loci of $\mathscr{M}_{\bar X}$ near $[\bar X]$, we recall Friedman's result on (partial) smoothing of ODPs: 

\begin{proposition} \cite{rF}
Let $w_i = a_{i1} r_1 + \ldots + a_{i\mu} r_\mu$, then the divisor $D_i := \{w_i = 0\} \subset \M_{\bar X}$ is the loci where the sphere $S_i$ shrinks to an ODP $p_i$. 
\end{proposition}

It is clear that the discriminant loci $D_B = \bigcup_{i = 1}^k D_i$ is not a normal crossing divisor. Rather it is a \emph{central hyperplane arrangement}. 

Under a suitable choice of homology symplectic basis, the $\beta$-periods in the transversal directions are given by 
$$
u_p = \p_p u = \int_{\beta_p} \Omega
$$ 
for some function $u$. The Bryant--Griffiths--Yukawa couplings are then extended over the boundary $D_B$ and satisfy
$$
u_{pmn} := \p^3_{pmn} u = O(1) + \sum_{i = 1}^k \frac{1}{2\pi \sqrt{-1}} \frac{a_{ip}a_{im}a_{in}}{w_i}
$$
for $1 \le p, m, n \le \mu$. It is holomorphic if one of the indices is outside this range.

The collection $\{u_{pmn}\}$ is the essential part of the Gauss--Manin connection $\nabla^{GM}$ on $\M_X$ which has regular singular extension over $D_B$.

Similarly, let $u = \sum_{p = 1}^\rho u^p T_p \in H^2(Y)/H^2(X)$, $D^i := \{ \sum\nolimits_{p = 1}^\rho b_{ip} u^p = 0 \}$, $i = 1, \ldots, k$. By the multiple cover formula of GW invariants we know that $QH(Y)$, or its Dubrovin connection, is regular singular along $D^A = \bigcup D^i$.

Let $y = \sum_{i = 1}^k y_i e_i \in \Bbb C^k$, with $e^1, \ldots, e^k$ being the dual basis on $(\mathbb{C}^k)^\vee$. The trivial logarithmic connection on $\underline{\mathbb{C}}^k \oplus (\underline{\mathbb{C}}^k)^\vee \longrightarrow \mathbb{C}^k$ is defined by
\begin{equation*} 
\nabla^k = d + \frac{1}{z} \sum_{i = 1}^k \frac{d y_i}{y_i} \otimes (e^i \otimes e_i^*).
\end{equation*}
The statement $A^t B = 0$ in Theorem \ref{t:bes} leads to an orthogonal sum
\begin{equation} \label{e:AB}
\Bbb C^k = {\rm image}\,A \mathop{\oplus}^\perp {\rm image}\,B \cong V^* \oplus H^2(Y)/H^2(X).
\end{equation}

\begin{theorem} \cite[Theorem 4.1]{LLW3} \label{t:trivial} {\ }
Under the identification \eqref{e:AB},
\begin{itemize}
\item[(1)] when restricted to $V^*$, $\nabla^k$ is naturally identified with the logarithmic (regular singular) part of $\nabla^{GM}$;

\item[(2)] when restricted to $H^2(Y)/H^2(X)$, $\nabla^k$ is naturally identified with to the logarithmic part of $\nabla^{\rm Dubrovin}$.
\end{itemize}
\end{theorem}

\subsection{Global aspects}

Denote by $\mathscr{A}(-)$ the GW theory and $\mathscr{B}(-)$ the variations of Hodge structure. Theorem \ref{t:trivial} provides evidence to
\begin{equation*}
\mbox{``excess $\mathscr{A}$ theory'' + ``excess $\mathscr{B}$ theory'' = ``trivial''}
\end{equation*}
through the partial exchange of quantum information attached to \emph{vanishing cycles} on both the $\mathscr{A}$ and $\mathscr{B}$ theories. For the full information on quantum $\mathscr{A}$, $\mathscr{B}$ theories, we proved the following result:

\begin{theorem} \cite[Theorem 0.3]{LLW3}
Let $[X]$ be a nearby point of $[\bar{X}]$ in $\M_{\bar{X}}$.
\begin{itemize}
\item[(1)] 
The theory $\mathscr{A}(X)$ is a sub-theory of $\mathscr{A}(Y)$ (e.g.~quantum sub-ring in genus 0).
\item[(2)]
The theory $\mathscr{B}(Y)$ is a sub-theory of $\mathscr{B}(X)$ (invariant sub-VHS).
\item[(3)]
The theory $\mathscr{A}(Y)$ can be reconstructed from a ``refined $\mathscr{A}$ theory'' on 
$$
X^{\circ} := X \setminus \bigcup\nolimits_{i=1}^k S_i
$$ 
``linked'' by the vanishing spheres in $\mathscr{B}(X)$.
\item[(4)] 
The theory $\mathscr{B}(X)$ can be reconstructed from the VMHS on $H^3(Y^\circ)$, 
$$
Y^{\circ} := Y \setminus \bigcup\nolimits_{i=1}^k C_i,
$$ 
``linked'' by the exceptional curves in $\mathscr{A}(Y)$.
\end{itemize}
\end{theorem}

The definition of the linked GW invariant in (3) is really a reformulation of the \emph{discreteness of components} appearing in the \emph{virtual cycle} form of the degeneration formula for conifold transitions of Calab--Yau 3-folds:
$$
\langle -\rangle_{g, \beta}^X = \sum\nolimits_{\gamma \mapsto \beta} \langle -\rangle_{g, \gamma}^Y.
$$
The sum is a \emph{finite sum}. However, no method is known to single out the individual term in it. To get $QH(Y)$ from $QH(X)$, it requires a blow-up formula of GW invariants where the blow-up center is a Weil divisor.

Weil divisors on $\bar X$ can be constructed from the relation matrix $B$ on $\mathfrak{X}_t$. Indeed, $\sum_{i = 1}^k b_{ij}[S_i] = 0$ implies that there are real 4-chains $W_{j, t}$ such that 
$$
\sum\nolimits_{i = 1}^k b_{ij} S_i = \p W_{j, t}.
$$
When $t \to 0$, we get homology cycles $W_j := W_{j, 0}$, $j = 1, \ldots, \rho$, since now $\p W_j$ is supported at the ODPs. From our definition of Calabi--Yau varieties it can be shown that $W_j$'s are represented by algebraic cycles, hence they give rise to Weil divisors on $\bar X$. 

The projective small resolution $\psi: Y \to \bar X$ transforms all these $W_j$'s into Cartier divisors. Indeed, let $W$ be the sum of supports of all $W_j$'s. Then $Y = {\rm Bl}_W \bar X$, i.e.~the blow-up of the ideal sheaf $\mathscr{I}_W \subset \mathscr{O}_{\bar X}$. A non-Cartier Weil divisor is simply a \emph{non-complete intersection divisor}. Thus the problem is essentially a problem on finding \emph{a blow-up formula of GW theory with non-complete intersection center}. Notice that the GW theory on $\bar X$ is so far undefined in the literature. However from the deformation invariance of (log) GW theory one may in practice identify it with the GW theory on $X$. 

\begin{example} [Determinantal transitions \cite{LLW4}]
Let $Y \subset S \times P^n$ be the zero loci of sections $s_i \in \Gamma(S \times P^n, \mathscr{L}_i)$ where $\mathscr{L}_i \to S \times P^n$ are line bundles of the form $\mathscr{L}_i = L_i \boxtimes \mathscr{O}_{P^n}(1)$ with $L_i$ being semi-ample on $S$.    

Let $[x_0: \cdots : x_n]$ be the homogeneous coordinates on $P^n$. We write
\begin{equation} \label{e:sij}
s_i = \sum\nolimits_{j = 0}^n s_{ij}\, x_j, \qquad i = 0, \ldots, n,
\end{equation}
where $s_{ij} \in \Gamma(S, L_i)$. We are interested in studying the restriction of the projection map $\pi: S \times P^n \to S$ to $Y$. Define $\bar X = \pi(Y) \subset S$ and 
$$
\psi = \pi|_Y: Y \to \bar X.
$$ 
The variety $\bar X$ has defining equation 
$$
\Delta := \det s_{ij} = 0.
$$  
For  $p \in \bar X$, since $s_{ij}(p)$'s are fixed, $\psi^{-1}(p)$ is not unique if and only if equation \eqref{e:sij} has more than one dimensional solutions in $P^n$, i.e.~$p$ is a singular point of $\bar X$. The  contraction $\psi: Y \to \bar X$ is called a determinantal contraction. Notice that 
$$
\Delta \in \Gamma(S, \bigotimes\nolimits_{i = 0}^n L_i).
$$
If for general sections $\tau \in \Gamma(S, \bigotimes_{i = 0}^n L_i)$ the variety $X_\tau$ defined by $\tau = 0$ is smooth, then it gives rise to a transition $Y \searrow X$. If furthermore $\bar X$ has only ODPs, then we get a conifold transition. These properties hold for CICY 3-folds transitions which have been studied extensively in the literature. 

Given a determinantal transition, we proceed to determine the GW theory on $Y$ in terms of the one on $X$ and the data $L_i$'s. Our goal is to replace the \emph{extrinsic data} $L_i$'s by the \emph{intrinsic data} associated to $\psi$, namely the Weil divisor $W$ which gives $Y = {\rm Bl}_W \bar X$. Here we consider the $g = 0$ case:

Let $h = c_1(\mathscr{O}_{P^n}(1))$. As in \eqref{e:QL+QLH}, the quantum Lefschetz gives \begin{equation*}
\begin{split}
I^Y &= J^S \frac{\prod_{i = 0}^n (L_i + h)}{(h)^{n + 1}},\\
I^X &= J^S (\sum\nolimits_{i = 0}^n L_i),
\end{split}
\end{equation*}
where we omit the curve class $\beta$ in the subscript. Then
\begin{equation} \label{rel-IYX}
I^{Y/X} = \frac{\prod_{i = 0}^n (L_i + h)}{(h)^{n + 1} (\sum_{i = 0}^n L_i)}.
\end{equation}

We note that the divisor $h$ on $S \times P^n$ coming from $\mathscr{O}_{P^n}(1)$ restricts to a divisor, still called $h$, on $Y$. If $p \in \bar X$ is a point with positive dimensional fiber $\psi^{-1}(p) \subset \{p\} \times P^n$ then $h$ intersects $\psi^{-1}(p)$ non-trivially since $h$ comes from a hyperplane in $P^n$. When $\psi$ is a small contraction, this effective divisor $W = \psi_*(h) \subset \bar X$ is thus the Weil divisor we are seeking for, and then $Y$ is the blow-up of $\bar X$ along $W$. 

It remains to interpret the \emph{$K$-trivial factor} \eqref{rel-IYX} in term of the linear system $|W|$. As in the case of smooth blow-ups along complete intersection centers (cf.~Remark \ref{r:L-system}), we will give intrinsic meaning of $I^{Y/X}$ in terms of decompositions of linear systems. The details will appear in \cite{LLW4}.
\end{example}

\end{document}